# ON HEEGNER POINTS OF LARGE CONDUCTORS

C. KHARE AND C. S. RAJAN

ABSTRACT. Given a parametrisation of an elliptic curve by a Shimura curve, we show that the images of almost all Heegner points are of infinite order.

## 1. Preliminaries

Let $\mathcal{B}$ be an indefinite quaternion algebra over $\mathbf{Q}$ that is ramified at an even number of primes $l_1, \ldots, l_n$. Denote by $D = \prod_i l_i$ the reduced disriminant of $B$. Let $N \geq 1$ be an integer that is relatively prime to $D$. Let $\mathcal{O}(N)$ denote a $\mathbf{Z}$-Eichler order in $\mathcal{B}$ of level $N$, and let $X_0^D(N)$ denote the corresponding Shimura curve defined over $\mathbf{Q}$. We recall the modular interpretation of $X_0^D(N)(\mathbf{C})$ [BD, Section 4, page 473]. The points $X_0^D(N)(\mathbf{C})$ correspond to triples $(A, \iota, C)$ where $A$ is an abelian surface over $\mathbf{C}$, $\iota : \mathcal{O}(N) \hookrightarrow \mathrm{End}(A)$ is an embedding (quaternionic multiplication), and $C$ is a level $N$ structure, i.e., a subgoup of $A(\mathbf{C})$ isomorphic to $\mathbf{Z}/N\mathbf{Z}$ and cyclic under the action of $\mathcal{O}(N)$.

Let $K$ be an imaginary quadratic field of discriminant $\mathrm{disc}(K)$, such that $ND$ and $\mathrm{disc}(K)$ are coprime. We further assume that the primes dividing $N$ are split in $K$ and the primes dividing $D$ remain inert in $K$. Let $f$ be an integer prime to $ND$, and let $\mathcal{O}_{K,f}$ be the order in $K$ defined by $\mathbf{Z} + f\mathcal{O}_K$.

**Definition 1.** [BD, Definition 5.1, page 482] With the above notations, a Heegner point $P_{K,f}$ on $X := X_0^D(N)$ attached to the pair $(K, f)$ corresponds to a triple $(A, \iota, C)$, such that $\mathrm{End}_{\mathcal{O}(N)}(A) = \mathcal{O}_{K,f}$, where $\mathrm{End}_{\mathcal{O}(N)}(A)$ denotes the endomorphisms of $A$ which commute with the action of $\mathcal{O}(N)$ and preserves the level $N$ structure.

We have the ring class field $K_f$ of conductor $f$ with $\mathrm{Gal}(K_f/K) \simeq \mathrm{Pic}(\mathcal{O}_{K,f})$ such that $\mathrm{Gal}(K_f/\mathbf{Q})$ is dihedral, i.e., the action of the non-trivial element of $\mathrm{Gal}(K/\mathbf{Q})$ on $\mathrm{Gal}(K_f/K)$ is by inversion. Note that $K_f \cap \mathbf{Q}(\mu_\infty) = K$, where $\mathbf{Q}(\mu_\infty)$ is the maximal abelian extension of $\mathbf{Q}$ generated by all roots of unity. Thus these extensions $K_f/\mathbf{Q}$ are anti-cyclotomic. The fundamental property of Heegner points that is needed here is as follows.

**Proposition 1.** *The minimal field of definition of a Heegner point $P_{K,f}$ over $K$ is $K_f$.*

*Proof.* This follows from the theory of complex multiplication by Shimura and Taniyama [ST]. A proof in this specific context can also be deduced from Lemma 2.5 and Theorem 5.3 of [BD]. □

**Corollary 1.** *Given $m$ a natural number, there are only finitely many pairs $(K, f)$ such that the degree of the field of definition of $P_{K,f}$ over $K$ is at most $m$.*

*Proof.* We have the following exact sequence:
$$(\mathbf{Z}/f\mathbf{Z})^* \times W \to (\mathcal{O}_K/f\mathcal{O}_K)^* \to \mathrm{Pic}(\mathcal{O}_{K,f}) \to \mathrm{Pic}(\mathcal{O}_K) \to 1,$$
where $W$ is the group of roots of unity in $K$. By Brauer-Siegel theorem $|\mathrm{Pic}(\mathcal{O}_K)| \leq m$ for only finitely many imaginary quadratic fields $K$, and from the above seqeunce we deduce that $|\mathrm{Pic}(\mathcal{O}_{K,f})| \leq m$ for only finitely many pairs $(K, f)$. □

## 1.1. Modular parametrisations.

Let $E/\mathbf{Q}$ be an elliptic curve of conductor $ND$: thus by recent results (Wiles, Taylor-Wiles, Diamond, Conrad-Diamond-Taylor, Breuil-Conrad-Diamond-Taylor) there is a non-constant morphism $X_0(ND) \to E$ defined over $\mathbf{Q}$. From this by standard arguments of which the key component is the Jacquet-Langlands correspondence we deduce the existence of a non-constant morphism
$$\phi : X_0^D(N) \to E$$
defined over $\mathbf{Q}$.

## 1.2. Theorem.

We state the main theorem of the note.

**Theorem 1.** *$\phi(P_{K,f})$ is a non-torsion point of $E(K_f)$ for almost all pairs $(K, f)$.*

*Remark.* The analogous theorem in the setting of modular curves (i.e., for $D = 1$) was proved in [NS]. For this reason we exclude $D = 1$ from the arguments, although our method works in that case too. The method of [NS] on the other hand does not generalise to the case $D > 1$, as it uses the existence of cusps on the (classical) modular curves.

*Remark.* The methods employed in this paper should allow an extension of the main theorem to any abelian variety quotient of the Jacobian of the modular curve, by invoking Ribet's theorem on the image of Galois for modular abelian varieties and by a suitable extension of Proposition 4 for CM fields.

*Remark.* By further refining the arguments in this paper, it can be shown that the subgroup generated by $\phi(P_{K,f})$ in $E(\overline{\mathbf{Q}})$ as either $K$ or $f$ vary, is of infinite rank.

## 2. Proof of the theorem

Denote by $\mathbf{Q}(\mu_\infty)$ the maximal abelian extension of $\mathbf{Q}$. For any field $F$, let $E(F)_t$ denote the torsion subgroup of $E(F)$. Let $E_t$ be the torsion points of $E(\overline{\mathbf{Q}})$. For any extension $L$ of $\mathbf{Q}$ let $L(E_t)$ be the fixed field of the kernel of the action of $\mathrm{Gal}(\bar{\mathbf{Q}}/L)$ on $E_t$. Denote by $\Gamma_L$ the Galois group $\mathrm{Gal}(L(E_t)/L)$. There is a natural faithful action of $\mathrm{Gal}(L(E_t)/L)$ on $E_t$, and we denote by $\Gamma_L$ the image inside $GL_2(\hat{\mathbf{Z}})$ of $\mathrm{Gal}(L(E_t)/L)$. By the Weil pairing we have the inclusion $L\mathbf{Q}(\mu_\infty) \subset L(E_t)$.

### 2.1. $E$ does not have CM.

We have the following result of Serre [Se].

**Proposition 2.** *Let $E$ be an elliptic curve over $\mathbf{Q}$ without CM. Then $\Gamma_{\mathbf{Q}(\mu_\infty)} = \mathrm{Gal}(\mathbf{Q}(E_t)/\mathbf{Q}(\mu_\infty))$ is isomorphic to an open subgroup of $SL_2(\hat{\mathbf{Z}})$.*

**Corollary 2.** *The maximal abelian extension of $\mathbf{Q}(\mu_\infty)$ contained inside $\mathbf{Q}(E_t)$ is of finite degree over $\mathbf{Q}(\mu_\infty)$.*

**Proposition 3.** *Let $\mathcal{K}$ denote the compositum of the fields $K_f$, as $K$ varies over imaginary quadratic fields and $f$ varies over natural numbers. Then $\mathcal{K} \cap \mathbf{Q}(E_t)$ is a finite extension of $\mathbf{Q}$.*

*Proof.* The proof follows from the following facts:
1. $K \subset \mathbf{Q}(\mu_\infty)$.
2. Hence the compositum $L$ of $\mathbf{Q}(\mu_\infty)$ with $K_f \cap K(E_t)$, where $K$ runs over all imaginary quadratic fields and $f$ all natural numbers, is an abelian extension of $\mathbf{Q}(\mu_\infty)$ contained in $\mathbf{Q}(E_t)$.
3. By the corollary above $L/\mathbf{Q}(\mu_\infty)$ is a finite extension.
4. Since $K_f$ is anticyclotomic, $K_f \cap \mathbf{Q}(\mu_\infty) = K$.

□

**Corollary 3.** *There are only a finite number of pairs $(K, f)$ as above such that the degree of $K_f \mathbf{Q}(E_t)/\mathbf{Q}(E_t)$ is less than a given number $m$.*

*Proof.* This follows from the above proposition and Corollary 1, as this degree is the degree of $K_f$ over $K_f \cap \mathbf{Q}(E_t)$. □

**Proof of theorem when $E$ does not have CM.** Suppose $\phi(P_{K,f})$ is a torsion point on $E$. Then by Proposition 1, we deduce that $K_f\mathbf{Q}(E_t)/\mathbf{Q}(E_t)$ is bounded by the degree of $\phi$. Hence the theorem follows from the above corollary in the case when $E$ does not have CM.

2.2. **$E$ has CM.** Let $E$ be an elliptic curve with CM by an order in an imaginary quadratic field $L$. There are two types of pairs $(K, f)$ to consider.

*Case I.* $(K, f)$ with $K = L$.

**Lemma 1.** *The torsion of $E(L_f)$ is bounded independently of $f$.*

*Proof.* This is proved in Section 2.2 of [NS]. Note that any prime of $L$ lying over an inert prime of $\mathbf{Q}$, splits in $L_f$. Reducing $E$ modulo two such primes, the torsion of $E(L_f)$ injects into a fixed finite group, and the lemma follows. □

From this Theorem 1 follows for pairs $(L, f)$ as $\phi^{-1}(E(L_f)_t)$ is a finite set as $f$ varies over natural numbers, since $\phi$ is a finite map.

*Case II.* $(K, f)$ with $K \neq L$. Let $L^{ac}$ be the maximal anti-cyclotomic extension of $L$ and $L^{ab}$ be the maximal abelian extension of $L$.

The following proposition is essentially proved in Section 3 of [E].

**Proposition 4.** *Let $L$ be an imaginary quadratic field. Then there are only a finite number of pairs $(K, f)$ with $K$ an imaginary quadratic field not equal to $L$ and $f$ a natural number, such that the degree of $K_f L^{ab}/L^{ab}$ is less than a given number $m$.*

*Proof.* It is proved in [E, Section 3], that given a natural number $m$, there are only a finite number of pairs $(K, f)$ with $K \neq L$, such that the degree of $K_f L^{ac}/L^{ac}$ is less than $m$. By the linear disjointness of $K_f$ and $\mathbf{Q}(\mu_\infty)$ over $K$, the proposition follows. □

Arguing as in the non-CM case, the above corollary implies Theorem 1, when $E$ has CM and we are in *Case II*.

**Acknowledgement:** We thank N. Schappacher for raising the question of generalising the result of [NS] to the situation of Shimura curves, and for his interest in this work. We also thank D. Prasad for his useful comments regarding this paper.

## References

[BD]  M. Bertolini and H. Darmon, Heegner points, $p$-adic $L$-functions, and the Cerednik-Drinfeld uniformization, *Invent. Math.*, **131** (1998) 453-491.

[E]  B. Edixhoven, Special points on the product of two modular curves. *Compositio Math.* **114** (1998), no. 3, 315–328.

[NS]  J. Nekovár and N. Schappacher, On the asymptotic behaviour of Heegner points, preprint.

[Se]  J.-P. Serre, Propriétes galoisiennes des points d'ordre fini des courbes elliptiques, Inv. Math., **15**, (1972), 259-331.

[ST]  G. Shimura and Y. Taniyama, Complex multiplication of abelian varieties and its applications to number theory, Publ. Math. Soc. Japan **6** (1961).

School of Mathematics, Tata Institute of Fundamental Research, Homi Bhabha Road, Bombay - 400 005, INDIA.
 *E-mail address*: `shekhar@math.tifr.res.in`

School of Mathematics, Tata Institute of Fundamental Research, Homi Bhabha Road, Bombay - 400 005, INDIA.
 *E-mail address*: `rajan@math.tifr.res.in`